\documentclass[12pt,draftcls,onecolumn]{IEEEtran}
\IEEEoverridecommandlockouts

\usepackage{amsmath,amssymb,amsfonts,epsf,epsfig,times}
\usepackage[all]{xy}
\usepackage{graphicx,color}
\usepackage{cite}
\usepackage{authblk}
\usepackage{amssymb}
\usepackage{subfigure}
\usepackage{multirow}
\usepackage{gensymb}
\usepackage{tikz} \usetikzlibrary{arrows,calc,decorations.pathreplacing}
\usepackage{fullpage}
\begin{document}
\title{\LARGE{Escape Regions of the Active Target Defense Differential Game}}
\author{Eloy Garcia\thanks{Corresponding Author: \ttfamily{elgarcia@infoscitex.com}. \\ Eloy Garcia is a contractor (Infoscitex Corp.) with the Control Science Center of Excellence, Air Force Research Laboratory, Wright-Patterson AFB, OH 45433. \\ David Casbeer is with the Control Science Center of Excellence, Air Force Research Laboratory, Wright-Patterson AFB, OH 45433. \\ Meir Pachter is with the Department of Electrical Engineering, Air Force Institute of Technology,	Wright-Patterson AFB, OH 45433.}, David W. Casbeer, and Meir Pachter}

\newtheorem{theorem}{Theorem}
\newtheorem{lemma}{Lemma}
\newtheorem{proposition}{Proposition}

\maketitle 

\begin{abstract}
The active target defense differential game is addressed in this paper. In this differential game an Attacker missile pursues a Target aircraft. The aircraft is however aided by a Defender missile launched by, say, the wingman, to intercept the Attacker before it reaches the Target aircraft. Thus, a team is formed by the Target and the Defender which cooperate to maximize the separation between the Target aircraft and the point where the Attacker missile is intercepted by the Defender missile, while the Attacker simultaneously tries to minimize said distance. This paper focuses on characterizing the set of coordinates such that if the Target's initial position belong to this set then its survival is guaranteed if both the Target and the Defender follow their optimal strategies. Such optimal strategies are presented in this paper as well.
\end{abstract}

\section{INTRODUCTION} \label{sec:intro}
In multi-agent pursuit-evasion problems one or more pursuers try to maneuver and reach a relatively small distance with respect to one or more evaders, which strive to escape the pursuers. 
This problem is usually posed as a dynamic game \cite{Ganebny12}, \cite{Huang11}, \cite{Pham10}. Thus, a dynamic Voronoi diagram has been used in problems with several pursuers in order to capture an evader within a bounded domain \cite{Huang11}, \cite{Bakolas10}. On the other hand, \cite{Sprinkle04} presented a receding-horizon approach that provides evasive maneuvers for an Unmanned Autonomous Vehicle (UAV) assuming a known model of the pursuer's input, state, and constraints. In \cite{EarlDandrea07}, a multi-agent scenario is considered where a number of pursuers are assigned to intercept a group of evaders and where the goals of the evaders are assumed to be known. 
Cooperation between two agents with the goal of evading a single pursuer has been addressed in \cite{Fuchs10} and \cite{Scott13}.

In this paper we consider a zero-sum three-agent pursuit-evasion differential game. A two-agent team is formed which consists of a Target ($T$) and a Defender ($D$) who cooperate; the Attacker ($A$) is the opposition. The goal of the Attacker is to capture the Target while the Target tries to evade the Attacker and avoid capture. The Target cooperates with the Defender which pursues and tries to intercept the Attacker before the latter captures the Target. Cooperation between the Target and the Defender is such that the Defender will capture the Attacker before the latter reaches the Target. Such a scenario of active target defense has been analyzed in the context of cooperative optimal control in \cite{Boyell76}, \cite{Boyell80}. Indeed, sensing capabilities of missiles and aircraft allow for implementation of complex pursuit and evasion strategies \cite{Zarchan97}, \cite{Siouris04},
 and more recent work has investigated different guidance laws for the agents $A$ and $D$. Thus, in \cite{ratnoo11} the authors addressed the case where the Defender implements Command to the Line of Sight (CLOS) guidance to pursue the Attacker which requires the Defender to have at least the same speed as the Attacker.
In \cite{rubinsky13} the end-game for the TAD scenario was analyzed based on the minimization/maximization of the Attacker/Target miss distance for a \textit{non-cooperative} Target/Defender pair. The authors develop linearization-based Attacker maneuvers in order to evade the Defender and continue pursuing the Target. 
A different guidance law for the Target-Attacker-Defender (TAD) scenario was given by Yamasaki \textit{et.al.} \cite{Yamasaki10}, \cite{Yamasaki13}. These authors investigated an interception method called Triangle Guidance (TG), where the objective is to command the defending missile to be on the line-of-sight between the attacking missile and the aircraft for all time, while the Target aircraft follows some predetermined trajectory.
The authors show, through simulations, that TG provides better performance in terms of Defender control effort than a number of variants of Proportional Navigation (PN) guidance laws, that is, when the Defender uses PN to pursue the Attacker instead of TG. 

The previous approaches constrain and limit the level of cooperation between the Target and the Defender by implementing Defender guidance laws without regard to the Target's trajectory.

Different types of cooperation have been recently proposed in \cite{Perelman11}, \cite{Rusnak05}, \cite{Rusnak11}, \cite{Ratnoo12}, \cite{Shima11}, \cite{Shaferman10}, \cite{Prokopov13} for the TAD scenario.
In \cite{Rusnak11} optimal policies (lateral acceleration for each agent including the Attacker) were provided for the case of an aggressive Defender, that is, the Defender has a definite maneuverability advantage. A linear quadratic optimization problem was posed where the Defender's control effort weight is driven to zero to increase its aggressiveness. The work \cite{Ratnoo12} provided a game theoretical analysis of the TAD problem using different guidance laws for both the Attacker and the Defender. The cooperative strategies in \cite{Shima11} allow for a maneuverability disadvantage for the Defender with respect to the Attacker and the results show that the optimal Target maneuver is either constant or arbitrary. Shaferman and Shima \cite{Shaferman10} implemented a Multiple Model Adaptive Estimator (MMAE) to identify the guidance law and parameters of the incoming missile and optimize a Defender strategy to minimize its control effort. In the recent paper \cite{Prokopov13} the authors analyze different types of cooperation assuming the Attacker is oblivious of the Defender and its guidance law is known. Two different one-way cooperation strategies were discussed: when the Defender acts independently, the Target knows its future behavior and cooperates with the Defender, and vice versa. Two-way cooperation where both Target and Defender communicate continuously to exchange their states and controls is also addressed, and it is shown to have a better performance than the other types of cooperation - as expected.

Our preliminary work \cite{Garcia14}, \cite{Garcia15} considered the cases when the Attacker implements typical guidance laws of Pure Pursuit (PP) and PN, respectively. In these papers, the Target-Defender team solves an \textit{optimal control} problem that returns the optimal strategy for the $T-D$ team so that $D$ intercepts the Attacker and at the same time the separation between Target and Attacker at the instant of interception of $A$ by $D$ is maximized. The cooperative optimal guidance approach was extended (\cite{Pachter14Allerton}, \cite{Garcia15ACC}, \cite{Garcia15JGCD}) to consider a differential game where also the Attacker missile solves an optimal control problem in order to minimize the final separation between itself and the Target. In this paper, we focus on characterizing the region of the reduced state space formed by the agents initial positions for which survival of the Target is guaranteed when both the Target and the Defender employ their optimal strategies. The optimal strategies for each one of the three agents participating in the active target defense differential game are provided in this paper as well.

The paper is organized as follows. Section \ref{sec:Problem} describes the engagement scenario. Section \ref{sec:analysis} presents optimal strategies for each one of the three participants in order to solve the differential game discussed in the paper. The Target escape region is characterized in Section \ref{sec:escape}. Examples are given in Section \ref{sec:Example} and concluding remarks are made in Section \ref{sec:concl}.

\section{PROBLEM STATEMENT} \label{sec:Problem}
\begin{figure}
	\begin{center}
		\includegraphics[width=8.4cm,height=7.1cm,trim=1.6cm 1.1cm 1.4cm .2cm]{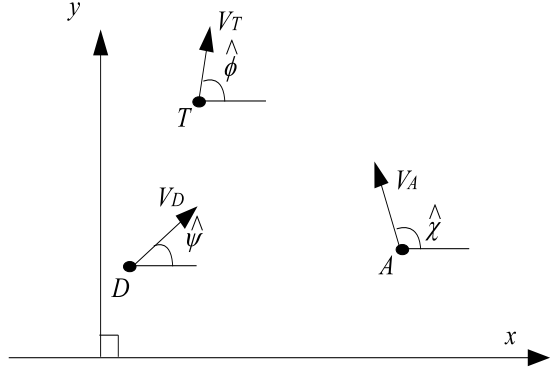}
	\caption{Active target defense differential game, $x_T>0$}
	\label{fig:problem description}
	\end{center}
\end{figure}

The active target defense engagement in the realistic plane $(x,y)$ is illustrated in Figure \ref{fig:problem description}. The speeds of the Target, Attacker, and Defender are denoted by $V_T$, $V_A$, and $V_D$, respectively, and are assumed to be constant. 
The simple-motion dynamics of the three vehicles in the realistic plane are given by:
\begin{align}
	\dot{x}_T&=V_T\cos\hat{\phi}, \ \ \ \ \ \ \ \ \ \ \ \dot{y}_T=V_T\sin\hat{\phi}    \label{eq:fixedT}  \\   
  \dot{x}_A&=V_A\cos\hat{\chi}, \ \ \ \ \ \ \ \ \ \ \  \dot{y}_A=V_A\sin\hat{\chi}   \label{eq:fixedADG} \\
	\dot{x}_D&=V_D\cos\hat{\psi}, \ \ \ \ \ \ \ \ \ \ \: \dot{y}_D=V_D\sin\hat{\psi}  	\label{eq:fixedD}
\end{align}
where the headings of $T$, $D$, and $A$ are, respectively, $\hat{\phi}$, $\hat{\psi}$, and $\hat{\chi}$.

In this game the Attacker pursues the Target and tries to capture it. The Target and the Defender cooperate in order for the Defender interpose himself between the Attacker and the Target and to intercept the Attacker before the latter captures the Target. Thus, the Target-Defender team searches for a cooperative optimal strategy, optimal headings $\hat{\phi}^*$ and $\hat{\psi}^*$, to maximize the separation between the Target and the Attacker at the time instant of the Defender-Attacker collision. The Attacker will devise its corresponding optimal strategy, optimal heading $\hat{\chi}^*$, in order to minimize the terminal $A-T$ separation/miss distance. Define the speed ratio problem parameter $\alpha=V_T/V_A$. We assume that the Attacker missile is faster than the Target aircraft, so that $\alpha <1$. In this work we also assume the Attacker and Defender missiles are somewhat similar, so $V_D=V_A$. 


In the following sections this problem is transformed to an aimpoint problem where each agent finds is optimal aimpoint. Furthermore, it is shown that the solution of the differential game involving three variables (the aimpoint of each one of the three agents) is equivalent to the solution of an optimization problem in only one variable.

\section{DIFFERENTIAL GAME} \label{sec:analysis}
We now undertake the analysis of the active target defense differential game. The Target ($T$), the Attacker ($A$), and the Defender ($D$) have ``simple motion" $\grave{\text{a}}$ la Isaacs \cite{Isaacs65}. We also emphasize that $T$, $A$, and $D$ have constant speeds of $V_T$, $V_A$, and $V_D$, respectively. We assume that $V_A=V_D$ and the speed ratio $\alpha=\frac{V_T}{V_A}<1$. We confine our attention to point capture, that is, the $D-A$ separation has to become zero in order for the Defender to intercept the Attacker.  $T$ and $D$ form a team to defend from $A$. Thus, $A$ strives to close in on $T$ while $T$ and $D$ maneuver such that $D$ intercepts $A$ before the latter reaches $T$ and the distance at interception time is maximized, while $A$ strives to minimize the separation between $T$ and $A$ at the instant of interception. Since the cost is a function only of the final time (the interception time instant) and the agents have simple motion dynamics, the optimal trajectories of each agent are straight lines.

\begin{figure}
	\begin{center}
		\includegraphics[width=8.4cm,height=7.1cm,trim=1.6cm 1.1cm 1.4cm .2cm]{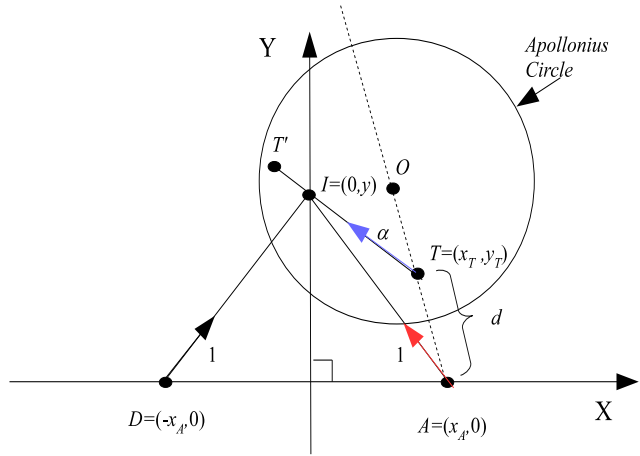}
	\caption{Active target defense differential game, $x_T>0$}
	\label{fig:oNE}
	\end{center}
\end{figure}

In Figure \ref{fig:oNE} the points $A$ and $D$ represent the initial positions of the Attacker and the Defender in the reduced state space, respectively. A Cartesian frame is attached to the points $A$ and $D$ in such a way that the extension to infinity of the segment $\overline{AD}$ in both directions represents the $X$-axis and the orthogonal bisector of $\overline{AD}$ represents the $Y$-axis. The state variables are $x_A$, $x_T$, and $y_T$. Notice that all points in the Left-Half-Plane (LHP) can be reached by the Defender before the Attacker does; similarly, all points in the Right-Half-Plane (RHP) can be reached by the Attacker before the Defender does.

In this paper we focus on the case where the Target is initially closer to the Attacker than to the Defender; in other words, assume that $x_T>0$. 

With respect to Figure \ref{fig:oNE} we note that the Defender will intercept the Attacker at point $I$ on the orthogonal bisector of $\overline{AD}$ at which time the Target will have reached point $T'$. The Attacker aims at minimizing the distance between the Target at the time instant when the Defender intercepts the Attacker, that is, the distance between point $T'$ and point $I$ on the orthogonal bisector of $\overline{AD}$ where the Defender intercepts the Attacker; the points $T$ and $T'$ represent the initial and terminal positions of the Target, respectively. 

\subsection{Cost/Payoff Function}  \label{subsec:xTg0cp}
\begin{figure}
	\begin{center}
		\includegraphics[width=8.4cm,height=6cm,trim=2.4cm 1.9cm 1.9cm .8cm]{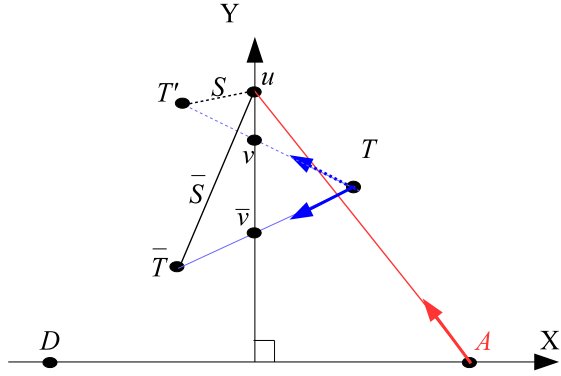}
	\caption{MaxMin Optimization Problem, $u>v$}
	\label{fig:MaxminEj}
	\end{center}
\end{figure}

When $x_T>0$ the Attacker and the Target are faced with a maxmin optimization problem: the Target chooses point $v$ and the Attacker chooses point $u$ on the Y-axis, see Figure \ref{fig:MaxminEj}. Additionally, the Defender tries to intercept the Attacker by choosing his aimpoint at point $w$ on the Y-axis.
Thus, the optimization problem is  
\begin{align}
    \max_{v,w} \ \min_{u} \ J(u,v,w).  
\end{align}
where the function $J(u,v,w)$ represents the distance between the Target terminal position and the point where the Attacker is intercepted by the Defender. The target tries to cross the orthogonal bisector of $\overline{AD}$ into the LHP where the Defender will be able to allow it to escape by intercepting the Attacker at the point $(0,u)$ on the orthogonal bisector of $\overline{AD}$. Therefore, the Defender's optimal policy is $w^*(u,v)=u$ in order to guarantee interception of the Attacker. The optimality of this choice by the Defender will be shown in Proposition \ref{prop:SPxtg0}.

\begin{figure}
	\begin{center}
		\includegraphics[width=8.4cm,height=6cm,trim=2.4cm 1.9cm 1.9cm .8cm]{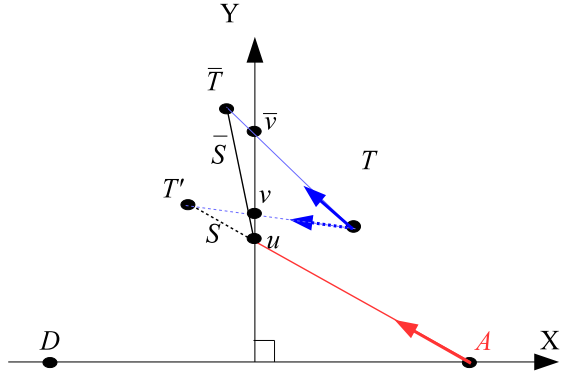}
	\caption{MaxMin Optimization Problem, $u<v$}
	\label{fig:MaxminEj2}
	\end{center}
\end{figure}

Since the Defender's optimal policy is $w^*=u$, the decision variables $u$ and $v$ jointly determine the distance $S$ between the Target terminal position $T'$ and the point $(0,u)$ where the Attacker is intercepted by the Defender. This distance is a function of the decision variables $u$ and $v$.
Thus, the Attacker and the Target solve the following optimization problem
\begin{align}
    \max_{v} \ \min_{u} \ J(u,v). 
\end{align}
Now, let us analyze the possible strategies. If the Target chooses $v$, the Attacker will respond and choose $u$. If $u\neq v$ the Target would correct his decision and choose some $\bar{v}$ such that $\bar{S}>S$, as shown in Figure \ref{fig:MaxminEj} for the case where $u>v$ and in Figure \ref{fig:MaxminEj2} for the case where $u<v$. In general, choosing $u\neq v$ is detrimental to the Attacker since his cost will increase. Thus, the Attacker should aim at the point $v$ which is chosen by the Target, that is, $u^*(v)=\arg \min_u J(u,v)=v$.   


\begin{proposition} Given the cost/payoff function $J(u,v)$, the solution $u^*$ and $v^*$ of the optimization problem $ \max_{v} \ \min_{u} J(u,v)$ is such that
\begin{align}
    u^*=v^*. 
\end{align}
Moreover, when $x_T>0$, the Attacker strategy is $u^*(v)=\arg \min_u J(u,v)=v$ so that it suffices to solve the optimization problem $\max_y J(y)$ where
\begin{align}
     J(y)=\alpha\sqrt{x_A^2+y^2}-\sqrt{(y-y_T)^2+x_T^2}.    \label{eq:Costy} 
\end{align}
\end{proposition}


\subsection{Critical Speed Ratio for Target Survival}
We assume that the Attacker is faster than the Target, for otherwise the Target could always escape without the help of the Defender. Thus, we assume that the speed ratio $0<\alpha<1$. Also, we assume that $x_T>0$. The Target needs to be able to break into the LHP before being intercepted by the Attacker for the Defender to be able to assist the Target to escape, by intercepting the Attacker who is on route to the Target. Thus, a solution to the active target defense differential game exists if and only if the Apollonius circle, which is based on the segment $\overline{AT}$ and the speed ratio $\alpha$, intersects the orthogonal bisector of $\overline{AD}$. This imposes a lower limit $\bar{\alpha}$ on the speed ratio, that is, we need $\bar{\alpha}<\alpha<1$. The critical speed ratio $\bar{\alpha}$ corresponds to the case where the Apollonius circle is tangent to the orthogonal bisector of $\overline{AD}$. And if the speed ratio $\alpha\geq 1$ the Target always escapes and there is no need for a Defender missile, that is, no target defense differential game is played out. The optimal strategies for the case $x_T<0$ can be obtained in a similar way as shown in this paper and the critical speed ratio is $\bar{\alpha}=0$.

\begin{proposition}
Assume that $x_T>0$. Then, the critical speed ratio $\bar{\alpha}$ is a function of the positions of the Target and the Attacker and is given by
\begin{align}
     \bar{\alpha}=\frac{\sqrt{(x_A+x_T)^2+y_T^2}-\sqrt{(x_A-x_T)^2+y_T^2}}{2x_A}.    \label{eq:alphasol}  
\end{align}
\end{proposition}
\textit{Proof}. The Attacker's initial position, the Target's initial position, and the center $O$ of the Apollonius circle are collinear and lie on the dotted straight line in Figure \ref{fig:oNE} whose equation is
\begin{align}
	  y=-\frac{y_T}{x_A-x_T}x + \frac{x_Ay_T}{x_A-x_T}. \nonumber
\end{align}
The geometry of the Apollonius circle is as follows: The center of the circle, denoted by $O$, is at a distance of $\frac{\alpha^2}{1-\alpha^2}d$ from $T$ and its radius is $\frac{\alpha}{1-\alpha^2}d$, where $d$ is the distance between $A$ and $T$ and is given by
\begin{align}
   d=\sqrt{(x_A-x_T)^2+y_T^2}. 
\end{align}
Hence, the following holds
\begin{align}
	   \Big(\frac{x_Ty_T}{x_A-x_T}-\frac{y_T}{x_A-x_T}x_0\Big)^2 + (x_0-x_T)^2  
		=\frac{\alpha^4}{(1-\alpha^2)^2}[(x_A-x_T)^2+y_T^2]   
\end{align}
and we calculate the coordinates of the center of the Apollonius circle
\begin{align}
   \left.
	 \begin{array}{l l}
     x_O=\frac{1}{1-\alpha^2}x_T-\frac{\alpha^2}{1-\alpha^2}x_A, \ \
		y_O=\frac{1}{1-\alpha^2}y_T.
	\end{array}  \label{eq:circCenter}  \right.
\end{align}
Consequently, the critical speed ratio $\bar{\alpha}$ is the positive solution of the quadratic equation 
\begin{align}
     x_T-\alpha^2x_A=\alpha\sqrt{(x_A-x_T)^2+y_T^2}    \label{eq:alphaeq}  
\end{align}
which is given by \eqref{eq:alphasol}.  
\begin{flushright}
$\square$
\end{flushright}

In general, it can be seen from Figure \ref{fig:oNE} that if $x_T<0$ then $\bar{\alpha}=0$ as well. 
We will assume $\bar{\alpha}<\alpha<1$, so that a solution to the active target defense differential game exists; otherwise, if $\alpha\leq\bar{\alpha}$, the Defender will not be able to help the Target by intercepting the Attacker before the latter inevitably captures the Target; and if $\alpha\geq 1$ then the Target can always evade the Attacker and there is no need for a Defender.

\subsection{Optimal Strategies}    \label{subsec:optimal}
When the Target is on the side of the Attacker, the Target chooses its aimpoint, denoted by $I$, on the orthogonal bisector of $\overline{AD}$ in order to maximize its payoff function \eqref{eq:Costy},
the final separation between Target and Attacker, and where $y$ represents the coordinate of the aimpoint $I$ on the orthogonal bisector of $\overline{AD}$. This is so because the Attacker will aim at the point $I$. In order to minimize \eqref{eq:Costy} the optimal strategy of the Attacker is to choose the same aimpoint $I$ on the orthogonal bisector of $\overline{AD}$, where it will be intercepted by the Defender. 

In order to find the maximum of \eqref{eq:Costy} we differentiate eq. \eqref{eq:Costy} in $y$ and set the resulting derivative equal to zero
\begin{align}
     \frac{dJ(y)}{dy}=\frac{\alpha y}{\sqrt{x_A^2+y^2}}-\frac{y-y_T}{\sqrt{(y-y_T)^2+x_T^2}}=0.    \label{eq:FderCosty}  
\end{align}
The following quartic equation in $y\geq 0$ is obtained
\begin{align}
     (1-\alpha^2)y^4 - 2(1-\alpha^2)y_Ty^3  
		+ \big((1\!-\!\alpha^2)y_T^2\!+\!x_A^2\!-\!\alpha^2x_T^2\big)y^2 
		- 2x_A^2y_Ty+x_A^2y_T^2=0.   \label{eq:Quartic}
\end{align}

In the sequel we focus on the case $0<\alpha<1$. In addition and without loss of generality assume that $y_T>0$.
Let us divide both sides of eq. \eqref{eq:Quartic} by $y_T^4$ and set $x_A=\frac{x_A}{y_T}$, $x_T=\frac{x_T}{y_T}$, and $y=\frac{y}{y_T}$, whereupon the quartic equation \eqref{eq:Quartic} assumes the canonical form 
\begin{align}
     (1-\alpha^2)y^4 - 2(1-\alpha^2)y^3 
		+ \big(1\!-\!\alpha^2 + x_A^2\!-\!\alpha^2x_T^2\big)y^2 
		- 2x_A^2y+x_A^2=0.
	  \label{eq:QuarticCan} 
\end{align}
We are interested in the real and positive solutions $y>0$ of the canonical quartic equation \eqref{eq:QuarticCan}. Eq. \eqref{eq:QuarticCan} has two real solutions,
\begin{align}
	 0<y_{R_1}<1 \ \  \textsl{and} \ \  y_{R_2}>1.  \nonumber
\end{align}
When $x_T=0$, \eqref{eq:QuarticCan} has two repeated solutions at $y=1$ and two complex solutions $y=\pm i\frac{1}{\sqrt{1-\alpha^2}}x_A$.

\textit{Remark}. Writing the quartic equation \eqref{eq:Quartic} as $f(y)=0$ we see that $f(0)=x_A^2y_T^2>0$, $f(y_T)=-\alpha^2 x_T^2y_T^2<0$, and $f(\infty)=+\infty$. Therefore, equation \eqref{eq:Quartic} has two real solutions. Equation \eqref{eq:Quartic} has a real solution $0<y<y_T$ and an additional real solution $y_T<y$, provided that $x_T\neq 0$. 
Note that the quartic equation \eqref{eq:Quartic} is parameterized by $x_T^2$, so whether $x_T>0$ or $x_T<0$ makes no difference as far as the solutions to the quartic equation \eqref{eq:Quartic} are concerned. However, if $x_T<0$ the applicable real solution is $y<y_T$, whereas if $x_T>0$ the applicable real solution is $y>y_T$. 

When $x_T>0$, by choosing his heading, the Target (and the Defender) thus choose the coordinate $y$ to maximize $J(y)$; that is, $y$ is the Target's (and Defender's) choice. Then the payoff is given by eq. \eqref{eq:Costy} and the expression for $\frac{dJ(y)}{dy}$ was shown in \eqref{eq:FderCosty}. The second derivative of the payoff function
\begin{align}
   \frac{d^2J(y)}{dy^2}=\frac{\alpha x_A^2}{(x_A^2+y^2)^{3/2}}-\frac{x_T^2}{\big((y-y_T)^2+x_T^2\big)^{3/2}}.  \label{eq:SderCosty}  
\end{align}
The Target is choosing $y$ to maximize the cost $J(y)$. Now, the Attacker reacts by heading towards the point $I$ on the orthogonal bisector of $\overline{AD}$ where, invariably, he will be intercepted by the Defender. Both the Target and the Attacker know that the three points $T,I,T'$ must be collinear. The defender will not allow the Attacker to cross the orthogonal bisector because then the Attacker will start to close in on the Target.

The optimal coordinate $y^*$ is the solution of the quartic equation \eqref{eq:Quartic} such that the second-order condition for a maximum holds on $\frac{d^2J(y)}{dy^2}<0$. In view of \eqref{eq:FderCosty} we know that 
\begin{align}
   \frac{1}{\sqrt{(y-y_T)^2+x_T^2}} = \alpha\frac{y}{y-y_T} \frac{1}{\sqrt{x_A^2+y^2}}.    \label{eq:opteq}  
\end{align}
and inserting \eqref{eq:opteq} into \eqref{eq:SderCosty} yields
\begin{align}
   \frac{d^2J(y)}{dy^2}=\frac{\alpha}{(x_A^2\!+\!y^2)^{3/2}} \Big(x_A^2 - \alpha^2\Big(\frac{y}{y\!-\!y_T}\Big)^3x_T^2\Big).   \label{eq:SderNeg}  
\end{align}
We have that $\frac{d^2J(y)}{dy^2}<0$ if and only if  
\begin{align}
   \frac{1}{\alpha^2}\Big(\frac{x_A}{x_T}\Big)^2 < \Big(\frac{y}{y\!-\!y_T}\Big)^3.   \label{eq:SderCond}  
\end{align}
Hence, the first real solution $y_1<y_T$ of the quartic equation \eqref{eq:Quartic} does not fulfill the role of yielding a maximum and the second real solution $y_2>y_T$ of \eqref{eq:Quartic} is the candidate solution. It is the Target who chooses $y^*$ to maximize the payoff $J(y)$. Note that $y_1=y_{R_1}y_T$ and $y_2=y_{R_2}y_T$, where $y_{R_1}$ and $y_{R_2}$ are the real solutions of \eqref{eq:QuarticCan}. 

\begin{figure}
	\begin{center}
		\includegraphics[width=8.4cm,height=6cm,trim=1.2cm 1.1cm 1.0cm .6cm]{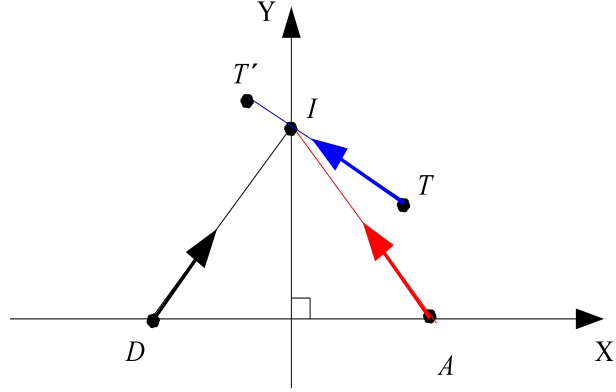}
	\caption{Optimal Play when $x_T>0$}
	\label{fig:OpSol}
	\end{center}
\end{figure}

Inserting eq. \eqref{eq:opteq} into eq. \eqref{eq:Costy} yields the Target and Defender payoff
\begin{align}
   \left.
	 \begin{array}{l l}
    J^*(y) &= \alpha\sqrt{x_A^2+y^2} - \frac{1}{\alpha}\frac{y-y_T}{y}\sqrt{x_A^2+y^2} \\
		 &=\frac{1}{\alpha}\sqrt{x_A^2+y^2}\big(\frac{y_T}{y} - (1-\alpha^2)\big) 
		\end{array}  \label{eq:payoff}  \right.
\end{align}
and using $y=y_2$
\begin{align}
   \left.
	 \begin{array}{l l}
    J^*(y)= \frac{1}{\alpha}\sqrt{x_A^2+y_2^2}\big(\frac{y_T}{y_2} - (1-\alpha^2)\big). 
	\end{array}  \label{eq:payoff2}  \right.
\end{align}
When $\alpha>\bar{\alpha}$ we have that $J^*(y)>0$. Hence, the solution $y_2>y_T$ of the quartic equation \eqref{eq:Quartic} must satisfy 
\begin{align}
   \left.
	 \begin{array}{l l}
   y_T<y_2<\frac{1}{1-\alpha^2} y_T.   \label{eq:OptCondy}  
	\end{array} \right.
\end{align}
This situation is illustrated in Figure \ref{fig:OpSol} where the three points $T$, $I$, and $T'$ are collinear. Concerning expression \eqref{eq:SderCond}, we also need the solution of the quartic equation \eqref{eq:Quartic} to satisfy
\begin{align} 
   y_2 < \frac{1}{1-\alpha^{2/3}(\frac{x_T}{x_A})^{2/3}} y_T.   \label{eq:OptCondy2}  
\end{align}
The second real solution $y_2$ of the quartic equation \eqref{eq:Quartic} must satisfy

\begin{align}
  y_T< y_2 <\min\left\{\frac{1}{1-\alpha^2} y_T, \frac{1}{1-\alpha^{2/3}(\frac{x_T}{x_A})^{2/3}} y_T\right\}.   \label{eq:OptCondymin}  
\end{align}
The points of intersection of the Apollonius circle with the $y$-axis (the orthogonal bisector) are $(0,\underline{y})$ and $(0,\overline{y})$, where $(0,\underline{y})$ and $(0,\overline{y})$ are the solutions of the quadratic equation
\begin{align}
   x_O^2+(y-y_O)^2=\frac{\alpha^2}{(1-\alpha^2)^2} d^2 \label{eq:yyQuad}	
\end{align}
where the distance $d=\sqrt{(x_A-x_T)^2+y_T^2}$ and the Apollonius circle's center coordinates are given by \eqref{eq:circCenter}. We have that
\begin{align}
   (y-y_O)^2=\frac{1}{(1-\alpha^2)^2}\big(\alpha^2d^2 - (x_T-\alpha^2x_A)^2\big) \nonumber	
\end{align}
where 
\begin{align}
\alpha^2d^2 - (x_T-\alpha^2x_A)^2>0  \nonumber
\end{align}
because $\alpha>\bar{\alpha}$ and, from \eqref{eq:alphaeq}, we have that $\bar{\alpha}d = x_T-\bar{\alpha}^2x_A$. Hence, 
\begin{align}
   y=y_O \pm \frac{1}{1-\alpha^2}\sqrt{\alpha^2d^2 - (x_T-\alpha^2x_A)^2} \nonumber	
\end{align}
which results in 
\begin{align}
   \underline{y}= \frac{1}{1-\alpha^2}\Big(y_T - \sqrt{\alpha^2y_T^2 + (1-\alpha^2)(\alpha^2x_A^2-x_T^2)}\Big)  \\
		\overline{y}= \frac{1}{1-\alpha^2}\Big(y_T + \sqrt{\alpha^2y_T^2 + (1-\alpha^2)(\alpha^2x_A^2-x_T^2)}\Big).	
\end{align}
The Target's choice of the optimal $y^*$, namely, the solution $y_2$ of the quartic equation \eqref{eq:Quartic} must satisfy the inequalities
\begin{multline}
	\frac{1}{1-\alpha^2}\Big(y_T \!-\! \sqrt{\alpha^2y_T^2 + (1\!-\!\alpha^2)(\alpha^2x_A^2-x_T^2)}\Big) 
	\leq y_2  \\
		\leq \frac{1}{1-\alpha^2}\Big(y_T \!+\! \sqrt{\alpha^2y_T^2 + (1\!-\!\alpha^2)(\alpha^2x_A^2-x_T^2)}\Big).	
\end{multline}

\begin{proposition} \label{prop:SPxtg0} (Saddle point equilibrium). Consider the case $x_T>0$. The strategy $y^*$ of the Target, where $y^*$ is the real solution of the quartic equation \eqref{eq:Quartic} which maximizes \eqref{eq:Costy}, and the strategy of the Defender of heading to the point $(0,y^*)$, together with the strategy of the Attacker of aiming at the point $(0,y^*)$, constitute a strategic saddle point,  that is 
\begin{align}
 \left.
	 \begin{array}{l l}
   &\left\{J(u^*,v^*,w),J(u^*,v,w^*),J(u^*,v,w)\right\} 
	  <J(u^*,v^*,w^*)<J(u,v^*,w^*) . \nonumber 
\end{array}   \right. 
\end{align}
\end{proposition}

\section{ESCAPE REGION}  \label{sec:escape}
In this section we analyze the Target's escape region for given Target and Attacker speeds, $V_T$ and $V_A$, respectively. In other words, for given speed ratio $\alpha=V_T/V_A$. Consider the active target defense differential game where the Attacker and the Defender missiles have the same speeds. When $x_T>0$, the critical value of the speed ratio parameter $\bar{\alpha}$ can be obtained as a function of the Attacker's and the Target's coordinates $x_A$, $x_T$, and $y_T$ such that the Target is guaranteed to escape since the Defender will be able to intercept the Attacker before the latter reaches the Target.

Now, for a given Target's speed ratio, $\alpha$, and for given Attacker's initial position, $x_A$, we wish to characterize the region of the reduced state space for which the Target is guaranteed to escape. In other words, we want to separate the reduced state space into two regions: $R_e$ and $R_{e_o}$. The region $R_e$ is defined as the set of all coordinate pairs $(x,y)$ such that if the Target's initial position $(x_T,y_T)$ is inside this region, then, it is guaranteed to escape the Attacker if both the Target and the Defender implement their corresponding optimal strategies. The region $R_{e_o}$, represents all other coordinate pairs $(x,y)$ in the reduced state space where the Target's escape is not guaranteed. 

\begin{proposition}
For given speed ratio $\alpha$ and for given Attacker's initial position $x_A$ in the reduced state space, the curve that divides the reduced state space into the two regions $R_e$ and $R_{e_o}$ is characterized by the right branch of the following hyperbola (that is, $x>0$)
\begin{align}
	  \frac{x^2}{\alpha^2x_A^2} -\frac{y^2}{(1-\alpha)^2x_A^2} = 1.     \label{eq:hb}
\end{align}
\end{proposition}
\textit{Proof}.
The requirement for the Target to escape being captured by the Attacker is that the Apollonius circle intersects the Y-axis. The radius of the Apollonius circle is 
\begin{align}
	  r=\frac{\alpha}{1-\alpha^2}\sqrt{(x_A-x_T)^2+y_T^2} \label{eq:r}
\end{align}
and the X-coordinate of its center is
\begin{align}
	  x_O=\frac{1}{1-\alpha^2}(x_T-\alpha^2x_A) \label{eq:xo}
\end{align}
If $x_T>\alpha^2x_A$, we need $r>x_O$ for the Defender to be of any help to the Target. Thus, $x_A>0$, $y_T\geq 0$, and $x_T>0$ must satisfy the condition
\begin{align}
	  \frac{1}{1-\alpha^2}(x_T-\alpha^2x_A) < \frac{\alpha}{1-\alpha^2}\sqrt{(x_A-x_T)^2+y_T^2} \label{eq:cond}
\end{align}
equivalently,
\begin{align}
	  x_T-\alpha^2x_A < \alpha\sqrt{(x_A-x_T)^2+y_T^2} \label{eq:cond2}
\end{align}
which is also equivalent to
\begin{align}
	  \frac{x_A^2}{(\frac{x_T}{\alpha})^2} +\frac{y_T^2}{\big(\frac{\sqrt{1-\alpha^2}}{\alpha}x_T\big)^2} > 1.     \label{eq:cond3}
\end{align}
When the `greater than' sign in inequality \eqref{eq:cond3} is changed to `equal' sign, the resulting equation defines the curve that divides the reduced state space into regions $R_e$ and $R_{e_o}$. Additionally, since the symmetric case $y_T<0$ can be treated in a similar way as the case $y_T>0$, we do not need to restrict $y_T$ to be greater than or equal to zero. Thus, the coordinate pairs $(x,y)$ such that
\begin{align}
	  \frac{x_A^2}{(\frac{x}{\alpha})^2} +\frac{y^2}{\big(\frac{\sqrt{1-\alpha^2}}{\alpha}x\big)^2} = 1     \nonumber
\end{align}
can be written in the hyperbola canonical form shown in \eqref{eq:hb}.
\begin{flushright}
$\square$
\end{flushright}

\textit{Remark}. Note that for a given speed ratio $\alpha$, the family of hyperbolas characterized by different values of $x_A>0$ shares the same center which is located at $C=(0,0)$, and the same asymptotes which are given by the lines $y=\frac{\sqrt{1-\alpha^2}}{\alpha}x$ and $y=-\frac{\sqrt{1-\alpha^2}}{\alpha}x$. One can also see that, for $0<\alpha<1$, the slope of the asymptotes increases as $\alpha$ decreases and viceversa. This behavior is expected since a relatively faster Target will be able to escape the Attacker when starting at the same position as a relatively slower Target.  

\begin{figure}
	\begin{center}
		\includegraphics[width=8.4cm,height=7.2cm,trim=2.4cm 1.5cm 4.4cm .2cm]{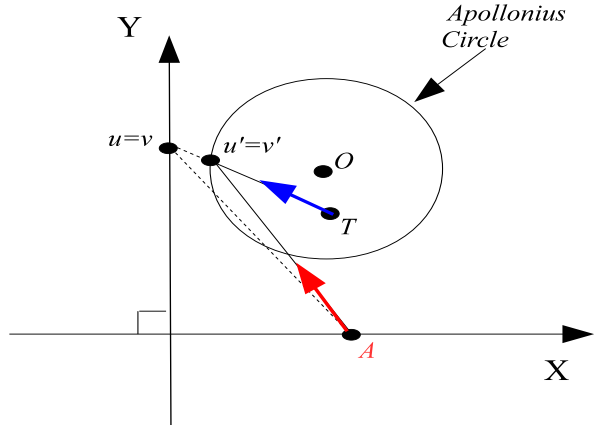}
	\caption{Attacker optimal strategy when $J(y^*)<0$}
	\label{fig:Jneg}
	\end{center}
\end{figure}

It is important to emphasize that if $(x_T,y_T)\in R_{e_o}$ then capture of the Target by the Attacker is guaranteed if the Attacker employs its optimal strategy. In this case the optimal strategies described in Section \ref{subsec:optimal} will result in $J(y^*)<0$. $J(y^*)$ being negative makes sense in terms of the differential game formulated in this paper (recall that the Attacker tries to minimize $J(y)$). However, the cost/payoff function $J(y)$ represents a distance and it does not make sense for it to be negative in the real scenario where the Attacker tries to capture the Target, i.e. the terminal separation $\overline{AT}$ should be zero instead of negative. Based on the solution of the differential game presented in Section \ref{subsec:optimal}, the Attacker is able to redefine its strategy and capture the Target, that is, to obtain $\overline{AT}=0$. The new strategy is as follows. The Attacker, by solving the differential game and obtaining the optimal cost/payoff, realizes that $J(y^*)<0$, then, it simply redefines its optimal strategy to be $u'^*(v')=v'$.

The details when the optimal strategies of Section \ref{subsec:optimal} result in $J(y^*)<0$ are as follows. $T$ chooses his aimpoint to be $v'$ that lies on the Apollonius circle. $A$ realizes that $J(y^*)<0$ (equivalently, the Apollonius circle does not intersect the Y-axis) and chooses his aimpoint $u'$ also on the Apollonius circle - see Figure \ref{fig:Jneg}. Given $T$'s choice of $v'$, the soonest $A$ can make $\overline{AT}=0$ is by capturing $T$ on the Apollonius circle (otherwise the Target will exit the Apollonius circle and the Defender may be able to assist the Target). Thus, $u'^*(v')=v'$. 

Similarly, $T$ solves the differential game and obtains $J(y^*)<0$. This information is useful to $T$ and it realizes that $D$ is unable to intercept $A$. Thus, $T$ will be prepared to apply passive countermeasures such as releasing chaff and flares. $T$ can also change its objective and find some $u'^*$ in order to optimize a different criterion such as to maximize capture time; however, this topic falls outside the scope of this paper.

\section{EXAMPLES} \label{sec:Example}

\begin{figure}
	\begin{center}
		\includegraphics[width=8.4cm,height=6.5cm,trim=.5cm .1cm .1cm .1cm]{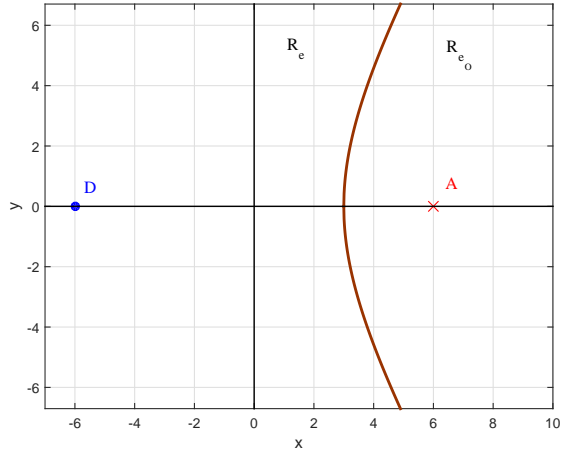}
	\caption{Example 1: Escape region}
	\label{fig:Reg1-1}
	\end{center}
\end{figure}

\begin{figure}
	\begin{center}
		\includegraphics[width=8.4cm,height=6.5cm,trim=.5cm .1cm .1cm .1cm]{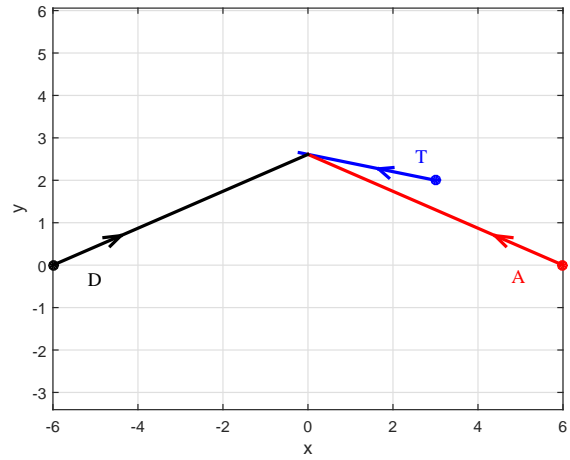}
	\caption{Example 1: Simulation}
	\label{fig:Ex1Sim}
	\end{center}
\end{figure}

\textit{Example 1}. Consider the speed ratio $\alpha=0.5$ and the Attacker's initial position $x_A=6$. The right branch hyperbola shown in Figure \ref{fig:Reg1-1} divides the Target escape/capture regions. 

Simulation: Let the Target's initial coordinates be $x_T=3$ and $y_T=2$. Note that $(x_T,y_T)\in R_e$. The Y-coordinate of the optimal interception point is given by $y^*=2.6108$. Figure \ref{fig:Ex1Sim} shows the results of the simulation. The optimal cost/payoff is $J(y^*)=0.2102$ and the Target escapes being captured by the Attacker. 

\begin{figure}
	\begin{center}
		\includegraphics[width=8.4cm,height=6.5cm,trim=.5cm .1cm .1cm .1cm]{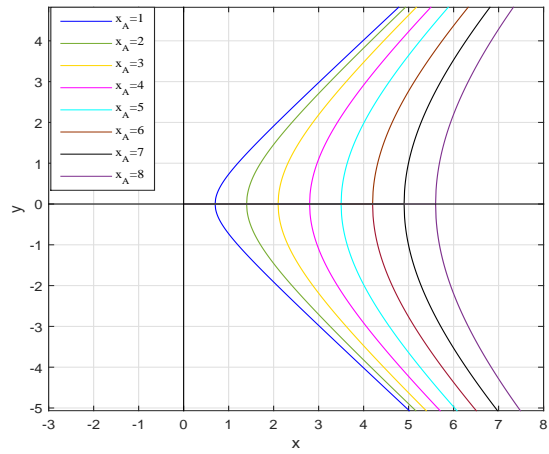}
	\caption{Example $2$: Family of right branch hyperbolas}
	\label{fig:RegM}
	\end{center}
\end{figure}

\textit{Example 2}. For a given speed ratio $\alpha$, we can plot a family of right hand hyperbolas on the same plane for different values of $x_A$. Consider $\alpha=0.7$. Figure \ref{fig:RegM} shows several hyperbolas \eqref{eq:hb} for values of $x_A=1,2,...,8$.

\section{CONCLUSIONS}  \label{sec:concl}
A cooperative missile problem involving three agents, the Target, the Attacker, and the Defender was studied in this paper. A differential game was analyzed where the Target and the Defender team up against the Attacker. The Attacker tries to pursue and capture the Target. The Target tries to evade the Attacker and the Defender helps the Target to evade by intercepting the Attacker before the latter reaches the Target. This paper provided optimal strategies for each one of the agents and also provided a further analysis of the Target escape regions for a given Target/Attacker speed ratio.


\end{document}